\documentclass [12pt] {report}
\usepackage{amssymb}
\newcommand {\D}[2] {\displaystyle\frac{\partial{#1}}{\partial{#2}}}

\newcommand {\Dd}[3] {\displaystyle\frac{\partial^2{#1}}{\partial{#2}\partial{#3}}}
\newcommand {\al} {\alpha}

\newcommand {\ga} {\gamma}
\newcommand {\la} {\lambda}

\newcommand {\de} {\delta}
\newcommand {\prtl} {\partial}
\newcommand {\fr} {\displaystyle\frac}
\newcommand {\wt} {\widetilde}
\newcommand {\wh} {\widehat}
\newcommand {\be} {\begin{equation}}
\newcommand {\ee} {\end{equation}}
\newcommand {\ba} {\begin{array}}
\newcommand {\ea} {\end{array}}
\newcommand {\bp} {\begin{picture}}
\newcommand {\ep} {\end{picture}}
\newcommand {\bc} {\begin{center}}
\newcommand {\ec} {\end{center}}
\newcommand {\bt} {\begin{tabular}}
\newcommand {\et} {\end{tabular}}
\newcommand {\lf} {\left}
\newcommand {\rg} {\right}

\newcommand {\cF} {{\cal F}}

\newcommand {\cI} {{\cal I}}

\newcommand {\cR} {{\cal R}}

\newcommand {\ses} {\medskip}

\newcommand {\e} {\mathop{\rm e}\nolimits}

\newcommand {\bibit} {\bibitem}
\newcommand {\nin} {\noindent}

\newcommand {\Ga} {\Gamma}
\newcommand {\De} {\Delta}

\usepackage{epsfig, latexsym,graphicx}
 \usepackage{amsmath}

\def\2#1#2#3{{#1}_{#2}\hspace{0pt}^{#3}}
\def\3#1#2#3#4{{#1}_{#2}\hspace{0pt}^{#3}\hspace{0pt}_{#4}}
\newcounter{sctn}
\def\sec#1.#2\par{\setcounter{sctn}{#1}\setcounter{equation}{0}
                  \noindent{\bf\boldmath#1.#2}\bigskip\par}

\begin {document}

\begin {titlepage}

\vspace{0.1in}

\begin{center}
{\Large \bf   Finsleroid-Finsler  Space of Involutive Case}

\end{center}

\vspace{0.3in}

\begin{center}

\vspace{.15in} {\large G.S. Asanov\\} \vspace{.25in}
{\it Division of Theoretical Physics, Moscow State University\\
119992 Moscow, Russia\\
{\rm (}e-mail: asanov@newmail.ru{\rm )}} \vspace{.05in}

\end{center}

\begin{abstract}

The Finsleroid-Finsler space is constructed over an underlying Riemannian space by the help of
a scalar $g(x)$ and an  input 1-form  $b$ of  unit length.
Explicit form of the entailed tensors, as well as the  respective spray coefficients, is evaluated.
The involutive case means the framework in  which the characteristic scalar $g(x)$ may vary in the direction
assigned by  $b$, such that $dg=\mu b$ with a scalar $\mu(x)$.
 We show by required calculation
 that the involutive case realizes through the
$A$-special relation
 the picture that instead of the
  Landsberg condition    $\dot A_{ijk}=0$ we have the vanishing
 $\dot{\al}_{ijk}=0$  with  the normalized tensor
 $\al_{ijk}=A_{ijk}/||A||$.
Under the involutive condition, the derivative tensor
 $A_{i|j}$
and the  curvature tensor
$R^i{}_k$
have explicitly  been found,
assuming the input 1-form $b$ be parallel.

\bigskip

\nin

\end{abstract}

\end{titlepage}

\vskip 1cm

{

\vskip 1cm

\setcounter{sctn}{1} \setcounter{equation}{0}

\bc  {\bf 1. Introduction    and synopsis  of new conclusions}  \ec

\bigskip

Among various possible methods to specify the Finsler space,
raising forth the Landsberg condition  $\dot A_{ijk}=0$
occupies an important geometrical role (see [1-3]).
 In the Finsleroid--Finsler space, the  condition
can be realized in a simple and attractive
way [4-7]. At the same time, the condition requires the Finsleroid charge $g$ to be a constant.
How should we overcome the restriction?

At the first sight, in the Finsler geometry  the weak Landsberg condition   $\dot A_{i}=0$
is to be considered as being a next--step extension of the proper Landsberg condition  $\dot A_{ijk}=0$.
 However, in the Finsleroid--Finsler space both the conditions
 are tantamount  (because of  the particular representation (1.13) of the Cartan tensor $A_{ijk}$).

A scrupulous analysis performed has revealed a remarkable observation that
an attractive method to permit $g\ne const$ is to use the nullification condition
 $\dot{\al}_{ijk}=0$ with  the normalized Cartan tensor  $\al_{ijk}$ (which is defined by (1.18)).
Clearly, the condition is attained when the $A$-special relation (2.5) holds.
Remarkably, the relation occurs being reachable upon assuming that the scalar $g(x)$ reveals
the involutive behaviour:    $dg=\mu(x)  b$ (see (3.1)).

In Section 2 we indicate the interesting implications of the $A$-special condition,
including the observation that
the skew--part of the $hv$-curvature tensor is proportional to
 the indicatrix curvature tensor (according to    (2.15)).
 It is the part that enters the right-hand side of the
 covariant conservation law (2.20).

In Section 3 the involutive case is formulated, showing that under the $b$-parallel condition
(which reads $\nabla b=0$, where
$\nabla$ means the Riemannian covariant derivative operative in the associated Riemannian space)
the case entails the $A$-special relation.

In Conclusions, several important ideas motivated our approach are emphasized.

Clear explicit representations of the basic tensors involved  are obtained systematically
in Appendix A by means of direct calculation.

In Appendix B, we indicate the explicit form for the spray coefficients of the space under study.
They include the part $E^k$ which involves the gradient of $g(x)$
(see (B.20) and (B.21)).

Appendix C is devoted to evaluations in the  involutive case.
Two key tensors have been explicitly  evaluated,
namely,
$A_{i|j}$ (given by (C.19))
and
$ R^i{}_k$
(given by (C.29)), assuming that the 1-form $b$ is parallel
(such that  $\nabla b=0$).

{

We  deal with the Finsler space notion which is specified by the condition that the Finslerian metric function
$K(x,y)$ be of the functional dependence
\be
K(x,y) =\Phi \lf(g(x), b_i(x), a_{ij}(x),y\rg)
\ee
of the particular case given by the formulas (A.36)--(A.41) of Appendix A.
In (1.1), the argument set $\bigl(g(x), b_i(x), a_{ij}(x)\bigr)$ involves, respectively,
a  scalar, a covariant vector field, and a Riemannian metric tensor.

The Finsleroid--Finsler space can be constructed as follows.
Let $M$ be an $N$-dimensional
$C^{\infty}$
differentiable  manifold, $ T_xM$ denote the tangent space to $M$ at a point $x\in M$,
and $y\in T_xM\backslash 0$  mean tangent vectors.
Suppose we are given on $M$ a positive--definite Riemannian metric
$${\cal S}=S(x,y).
$$
 Denote by
$$
\cR_N=(M,{\cal S})
$$
the obtained $N$-dimensional Riemannian space.
Let us also assume that the manifold $M$ admits a non--vanishing 1-form
$$
b= b(x,y)
$$
and
\be
|| b||_{\text{Riemannian}}=1.
\ee
It is convenient to  use the variable
\be
q =\sqrt{S^2-b^2}.
\ee

 The space $\cR_N$ entering the above definition is called the {\it associated Riemannian  space}.
 With respect to  natural local coordinates in the space
$\cR_N$
we have the local representations
\be
a^{ij}b_ib_j=1
\ee
and
\be
 b=b_i(x)y^i,
 \ee
 together with
 \be
  S= \sqrt{a_{ij}(x)y^iy^j}.
\ee
The reciprocity  $a^{in}a_{nj}=\de^i{}_j$ is assumed, where $\de^i{}_j$ stands for the Kronecker symbol.
The covariant index of the vector $b_i$  will be raised by means of the Riemannian rule
$$
 b^i=a^{ij}b_j,
 $$
 which inverse reads
 $$
  b_i=a_{ij}b^j.
  $$
We
also  introduce the tensor
\be
r_{ij}(x)~:=a_{ij}(x)-b_i(x)b_j(x)
\ee
to have the representation
\be
q=\sqrt{r_{ij}(x)y^iy^j}.
\ee
From (1.4) and (1.7) it follows that
\ses
\be
 r_{ij}b^j=0
\ee

{

From the fundamental metric function function $K$,
we  explicitly calculate  distinguished Finslerian tensors,
including
the covariant tangent vector $\hat y=\{y_i\}$,
the  Finslerian metric tensor $\{g_{ij}\}$
together with the contravariant tensor $\{g^{ij}\}$ defined by the reciprocity conditions
$g_{ij}g^{jk}=\de^k_i$, and the  angular metric tensor
$\{h_{ij}\}$, by making  use of the following conventional  Finslerian  rules in succession:
\be
y_i :=\fr12\D{K^2}{y^i}, \qquad
g_{ij} :
=
\fr12\,
\fr{\prtl^2K^2}{\prtl y^i\prtl y^j}
=\fr{\prtl y_i}{\prtl y^j}, \qquad
h_{ij} := g_{ij}-l_il_j,
\ee
where
\be
l_i=g_{ij}l^j, \qquad l^j=\fr{y^j}K.
\ee

After that, we can elucidate the algebraic structure of the associated  Cartan tensor
\be
A_{ijk} := \fr K2\D{g_{ij}}{y^k},
\ee
which leads to
the  simple representation
\be
A_{ijk}=\fr1N\lf(h_{ij}A_k+h_{ik}A_j+h_{jk}A_i-\fr1{A_hA^h}A_iA_jA_k\rg)
\ee
(see (A.87))
with
\be
A_k=g^{ij}A_{ijk}
\ee
and
\be
A_hA^h=\fr{N^2}{4}g^2
\ee
(see (A.63)).

Owing to (1.15), the norm
\be
||A||=\sqrt{A^kA_k}
\ee
is equal to
\be
||A||=\fr N2|g(x)|.
\ee

It is convenient to construct the {\it normalized Cartan tensor}
\be
\al_{ijk}:=\fr1{||A||}A_{ijk}
\ee
and the vector
\be
\al_{k}:=\fr1{||A||}A_{k}
\ee
which length is 1:
\be
\al_h\al^h=1.
\ee
We have
\be
\al_{ijk}=\fr1N\lf(h_{ij}\al_k+h_{ik}\al_j+h_{jk}\al_i-\al_i\al_j\al_k\rg)
\ee
everywhere in the Finsleroid--Finsler space.

{

In our analysis, an important role is played by the tensor
\be
{\cal H}_{ij}=h_{ij}-\al_i\al_j,
\ee
which obviously possesses the nullification properties
\be
{\cal H}_{ij}y^j =0, \qquad   {\cal H}_{ij}A^j=0.
\ee
The curvature of  indicatrix is well--known to be described by the tensor
\be
\hat R_i{}^j{}_{mn} := \fr1{K^2}\Bigl(\3Ahjm\3Aihn-\3Ahjn\3Aihm\Bigr).
\ee
In the Finsleroid--Finsler space, the  tensor possesses
the representation
\be
K^2\hat R_{ijmn}
=
\fr1{N^2}  (A^kA_k)
\Bigl(
 h_{in}h_{mj}  -h_{im}h_{nj}
\Bigr)
\ee
(see (A.93)).

\ses

In the next section we shall set forth an interesting special condition.

{

                        \bigskip

\setcounter{sctn}{2}
\setcounter{equation}{0}

\bc
 {\bf 2. $A$-special   condition}
 \ec

\bigskip

By means of the over-dot we denote the action of the operator
${\scriptstyle |m}
l^m$,
 such that
\be
\dot A_i= A_{i|m}l^m, \quad \dot A_{ijk}=  A_{ijk|m}l^m,
\qquad
\dot \al_i= \al_{i|m}l^m, \quad \dot \al_{ijk}=  \al_{ijk|m}l^m,
\ee
with ${\scriptstyle |m}$ meaning the $h$-covariant derivative (see (B.9)).

Let us    set forth   the nullification
\be
 \dot \al_{ijk}=0.
 \ee
Whenever  the representation (1.21) is valid, the condition (2.2) is equivalent to
the vanishing
 \be
\dot{\al}_i=0
\ee
of the normalized vector (1.19).

Denote
\be
\ga_k=\fr1{2A^hA_h}\,(A^mA_m)_{|k},   \qquad
\ga=\fr1{2A^hA_h}\,(A^mA_m)_{|k}l^k.
\ee
\ses
Assume that the {\it $A$-special relation}
\be
A_{i|k} =\ga_k    A_i
+  \eta{\cal H}_{ik}
\ee
holds,  where $ \eta$ is a scalar. The relation (2.5) can  obviously be written as
\be
\al_{i|k} =     \eta{\cal H}_{ik}.
\ee
Since  ${\cal H}_{ik}y^k=0$, from (2.5)  we directly conclude  that
\be
\dot A_i=\ga A_i,
\ee
which is obviously {\it tantamount to} (2.3).
From the representation (1.13) of the tensor
$A_{ijk}$ we obtain
\be
A_{ijk|l} ={\ga}_l  A_{ijk}
+ \eta \fr1N
 \bigl(     {\cal H}_{ij}{\cal H}_{kl}+{\cal H}_{ik}{\cal H}_{jl}+{\cal H}_{jk}{\cal H}_{il}    \bigr),
\ee
\ses
which entails
\be
\dot A_{ijk} = \ga   A_{ijk}.
\ee
Also, (2.3) entails the nullification
\be
\dot{\cal H}_{jk}=0
\ee
 (consider the representation (1.22) of the tensor ${\cal H}_{jk}$ and take into account that
$$
h_{jk|l}=0
$$
 in any Finsler space),
 where
 $$
 \dot{\cal H}_{jk} = {\cal H}_{jk|m}l^m.
 $$

{

The $hv$-curvature tensor
\be
P_{jikl}~: =-\bigl(A_{ijl|k}-   A_{jkl|i}  +  A_{kil|j}\bigr)
+A_{ij}{}^u\dot A_{ukl}
-
A_{jk}{}^u\dot A_{uil}
+A_{ki}{}^u\dot A_{ujl}
\ee
(this representation is tantamount to the definition (3.4.11) on p. 56 of the book [2])
gets reduced upon substituting (2.8) and (2.9):
\be
P_{jikl}=
-2\eta \fr1N
 \bigl(     {\cal H}_{ij}{\cal H}_{kl}+{\cal H}_{ik}{\cal H}_{jl}+{\cal H}_{jk}{\cal H}_{il}    \bigr)
+
\ga   (A_{ij}{}^u A_{ukl}  -  A_{jk}{}^u A_{uil}  +A_{ki}{}^u A_{ujl}).
\ee
Let us consider the skew--part
\be
P^{[ji]}{}_{kl}~:= \fr12(P^{ji}{}_{kl} -P^{ij}{}_{kl}).
\ee
From (2.12) it follows that
\be
P_{[ji]kl} =
\ga  ( A_{ki}{}^u A_{ujl} -  A_{jk}{}^u A_{uil} ).
\ee
In view of  the representation (1.24) of the indicatrix curvature tensor $\hat R_i{}^j{}_{mn}$,
we may write (2.14) as
\be
P_{[ji]kl} =
\ga K^2 \hat R_{jikl}.
\ee

Thus the following assertion is valid.

\ses

\ses

THEOREM 2.1. {\it If the $A$-special relation}
(2.5) {\it  holds together with the representation}
(1.13)
{\it of the Cartan tensor, then the skew--part of the $hv$-curvature tensor is proportional to
 the indicatrix curvature tensor, according to
}
 (2.15).

\ses

\ses

{

Since
\be
g^{jl}\Bigl(R_j{}^i{}_{il|t}+R_j{}^i{}_{lt|i}+R_j{}^i{}_{ti|l}\Bigr)=
P^l{}^i{}_{iu}R^u{}_{lt}+P^l{}^i{}_{lu}R^u{}_{ti}+P^l{}^i{}_{tu}R^u{}_{il}
\ee
(see the formula (3.5.3) on p. 58 of the book [2]), we have
\be
g^{jl}\Bigl(R_j{}^i{}_{il|t}+R_j{}^i{}_{lt|i}+R_j{}^i{}_{ti|l}\Bigr)=
2P^{[li]}{}_{iu}   R^u{}_{lt}-P^{[li]}{}_{tu}R^u{}_{li}.
\ee
\ses
so that the covariant divergence of   the tensor
\be
\rho_{ij}~:=\fr12(R_i{}^m{}_{mj}+R^m{}_{ijm})-\fr12g_{ij}R^{mn}{}_{nm}
\ee
is given by
\be
\rho^i{}_{j|i}
=-P^{[lm]}{}_{mu}   R^u{}_{lj}  + \fr12 P^{[lm]}{}_{ju}R^u{}_{lm}
\ee
which can be written as
\be
\rho^i{}_{j|i}=J_j
\ee
with
\be
J_j
=P^{[lm]}{}_{ku}
\Biggl(- R^u{}_{lj}\de^k{}_m  + \fr12 R^u{}_{lm} \de^k{}_j
\Biggr).
\ee

\ses

Using (2.15) together with (1.25) entails
\be
J_j=\fr14  g^2  \ga    \Bigl( h^l{}_uh^m{}_k  -h^l{}_kh^m{}_u  \Bigr)
 \lf(- R^u{}_{lj}\de^k{}_m  + \fr12 R^u{}_{lm} \de^k{}_j\rg).
\ee

{

                     \bigskip


\ses\ses

\setcounter{sctn}{3}
\setcounter{equation}{0}

\bc
 {\bf 3. Finsleroid--Finsler space upon  involution}
 \ec

\bigskip

Let us set forth the {\it involution condition}
\be
g_i=\mu b_i, \qquad \mu=\mu(x),
\ee
where $g_i=\partial g/\partial x^i$,
and formulate the following definition.

\ses

\ses

 {\large  Definition}.  The arisen  space
\be
\cI\cF\cF^{PD}_g~:=\{\cF\cF^{PD}_g ~ \text{with} ~ g_i=\mu b_i, ~ \mu=\mu(x)\}
\ee
is called the
{\it involutive Finsleroid--Finsler space}.
The involved $\mu(x)$ is called the {\it involution scalar}.

\ses

\ses

In the  space (3.2), the quantities defined in
(2.4) become simply
\be
\ga_k=\fr1gg_k,  \qquad
\ga=\fr1gg_kl^k,
\ee
so that
 the  $A$-special relation (2.5) takes on the form
\be
A_{i|k} = \fr1g g_k    A_i
+  \eta{\cal H}_{ik}.
\ee

We say that the space  $\cF\cF^{PD}_g $ is {\it $b$-parallel}, if the 1-form $b$ is parallel
in the sense of the associated Riemannian space, that is, when
\be
\nabla_ib_j=0.
\ee

\ses

It proves that the following theorem is valid.

\ses

\ses

THEOREM 3.1.  In the  $b$-parallel involutive   space  $\cI\cF\cF^{PD}_g$
the $A$-special relation  (3.4)  holds.

\ses

\ses

See Appendix C, in which all the involved evaluations (which are not short)
have been  presented
and the representation (C.19) has been arrived at,
 which belongs to the type (3.4).
The $\eta$ entered the right-hand part of (3.4) can be written down from (C.19).

As a direct consequence of the above theorem,
\be
\{\nabla b=0~\text{and}~dg=\mu  b\} ~ \Longrightarrow ~ \dot \al_i=0.
\ee

\ses

From (2.4) and (3.1) we have
\be
\ga=\fr1K\fr{\mu}g b.
\ee
With this formula,
the current (2.22) takes on the explicit representation
\be
J_j=\fr1{4K}\mu g b   \Bigl( h^l{}_uh^m{}_k  -h^l{}_kh^m{}_u  \Bigr)
 \lf(- R^u{}_{lj}\de^k{}_m  + \fr12 R^u{}_{lm} \de^k{}_j\rg)
\ee
which is proportional to the involution scalar $\mu$.

{

       \bigskip

\setcounter{sctn}{4}
\setcounter{equation}{0}

\bc
 {\bf 4. Conclusions}
 \ec

\bigskip

The Finsleroid--Finsler space involves a characteristic scalar, $g(x)$, such that
the vanishing of the scalar
reduces the space to a Riemannian space.
Varying $g(x)$ entails varying the form of the Finsleroid.
The Landsberg  case of the    Finsleroid--Finsler space  implies strictly $g=const$,
as a direct consequence  of the formulas (1.13)--(1.15).
To set a liberty to the scalar $g(x)$, we must overcome the restrictive case.
It proves that a fruitful idea is to substitute the condition
$\dot \al_{ijk}=0$
with   the Landsberg  condition $\dot A_{ijk}=0$ proper.
Would one assume $||A||=const$,  one  observes that
$\dot \al_{ijk}=0$  implies $\dot A_{ijk}=0$.
In the Finsleroid--Finsler space under study,   $g\ne const$ implies $||A||\ne const$
(see (1.17)).

 The involution condition (3.1) can be written as
$dg=\mu b_i(x)dx^i$ which means geometrically that the scalar $g(x)$  varies
in the direction assigned by the vector   $b_i(x)$.

The obtained involutive curvature tensor
$ R^i{}_k$
(given by (C.29))
 is of the novel type, being created by the gradient of the Finsleroid charge
and constructed from the involutive spray coefficients  $E^k$.
The tensor is meaningful   even if the associated Riemannian space is flat.

It would be   appealing to develop in future the extensions which can
 go over the $b$-parallel case $\nabla b=0$.

We have examined the conservation law for the fundamental tensor $\rho_{ij}$, obtaining the result (2.20)--(2.22).
In the Landsberg case, the $hv$-curvature tensor $P_{ijkl}$ is well-known to be totally symmetric in all four
of its indices
(see p. 60  in [2]), such that the skew-part $P^{[lm]}{}_{ku}$, and whence the right-hand part of (2.20),
 vanishes. Under the $A$-special condition, however, the tensor
 is meaningful, being proportional to the indicatrix curvature tensor in accordance with
 (2.15), so that the current $J_j$ given by (2.21) is not the zero.

Various   Finslerian ideas of applications  (see [8-10])
can well be matched to the ($g\ne const$)-Finsleroid-Finsler space.

{

            \bigskip

\ses

\setcounter{equation}{0}

\bc
{ \bf  Appendix A:     Evaluation of quantities of the  space ${\mathbf\cF^{PD}_g } $}
\ec

\ses\ses


Below, we evaluate the key objects of the  space ${\mathbf\cF^{PD}_g } $
under the general setting when the Finsleroid charge $g$  may depend on $x$, so that $g=g(x)$.
The unit norm $||b||=1$ of the input 1-form $b$ is assumed.
Our treatment will be of {\it local} character. Any dimension $N\ge 2$ is admissible.
The Riemannian squared length $S^2=b^2+q^2$ underlines the Finslerian space under study.
\ses

\ses

It is appropriate  to use the variables
\be
u_i~:=a_{ij}y^j,
\qquad
v^i~:=y^i-bb^i, \qquad v_m~:=u_m-bb_m=r_{mn}y^n\equiv a_{mn}v^n,
\ee
where $r_{mn}=a_{mn}-b_mb_n$.
We obtain
the relations
\be
r_{ij}=\D{v_i}{y^j},
\ee
\ses
\be
u_iv^i=v_iy^i=q^2,
 \qquad
v_ib^i=v^ib_i=0,
\ee
\ses
\be
 r_{in}v^n=v_i, \qquad    v_kv^k=q^2,
\ee
and
\be
\D b{y^i}=b_i, \qquad \D q{y^i}=\fr{v_i}q.
\ee

In terms of the variable
\be
w=\fr qb,
\ee
we obtain
\be
\D{w}{y^i}=\fr{z_i}{b^2q}, \qquad z_i=bv_i-q^2b_i
\equiv bu_i-S^2b_i,
\ee
and
$$
y^iz_i=0,\qquad
b^iz_i=b^2-S^2,
$$
together with
$$
 a^{ij}z_iz_j=S^2(S^2-b^2)\equiv S^2b^2\la,
$$
where
\be
\la=   w^2\equiv   \fr1{b^2}(S^2-b^2).
\ee

{

We also introduce  the $\eta$--{\it tensor}   by
means of the components
\be
\eta_{ij}~:=r_{ij}-\fr1{q^2}v_iv_j,
 \qquad \eta^i{}_j~:=r^i{}_j-\fr1{q^2}v^iv_j,
 \qquad \eta^{ij}~:=r^{ij}-\fr1{q^2}v^iv^j.
 \ee
It follows directly that
\be
\eta^n{}_j=a^{nm}\eta_{mj}, \qquad  \eta^{ij}= a^{in}\eta_n^j,
\ee
\ses
\be
\eta_{ni}y^i=0,
\ee
\ses
\be
 \eta_{ij}b^j=0,
\qquad
 \eta_{ij}z^j=0,
\qquad
a^{ij}\eta_{ij}=N-2,
 \ee
 \ses
and
\be
\D{\lf(\fr1qv_k\rg)}{y^j}=\fr1q\eta_{kj}, \qquad
\D{\eta_{ij}}{y^k}=-\fr1{q^2}(v_i\eta_{jk}+v_j\eta_{ik}).
\ee

We shall also use the vector
\be
e_k~:=\fr b{q^2}v_k-b_k
\equiv
 \fr b{q^2}u_k-\fr{S^2}{q^2}b_k,
\ee
obtaining
\be
e_k=-q\D{\lf(\fr bq\rg)}{y^k},
\ee
\ses
\be
\D{e_k}{y^j}=\fr b{q^2}\eta_{kj}-\fr {1}{q^2}v_ke_j=\fr b{q^2}\eta_{kj}-\fr {1}{b}(e_k+b_k)e_j,
\ee
\ses
\be
e_ky^k=0,
\ee
\ses
\be
e_kb^k=-\fr1{w^2}\la, \qquad e^j\eta_{ij}=0,
\qquad
\D{\eta_{ij}}{y^k}=-\fr1{b}(e_i\eta_{jk}+b_i\eta_{jk} +e_j\eta_{ik}+b_j\eta_{ik}),
\ee
and
\be
z_k=q^2e_k,
\ee
together with
\be
w^4 a^{ij}e_ie_j=(1+w^2)\la
\equiv
\fr{S^2}{b^2}\la.
\ee

{

\ses

Using the generating  function $V=V(x,w)$ defined from the representation
\be
K=bV,
\ee
we obtain
\be
\D K{y^i}=b_iV+\fr1{bq}z_iV', \qquad
\Dd{ K}{y^i}{y^j}=
\fr1q\eta_{ij}
V'
+\fr1{b^3q^2}z_iz_jV''.
\ee
The prime $\{'\}$ means differentiation with respect to $w$.
Taking into account the Finslerian rules
\be
l_i=\D K{y^i}, \qquad y_i=Kl_i,\qquad h_{ij}=K\Dd{ K}{y^i}{y^j}, \qquad
g_{ij}=h_{ij}+l_il_j,
\ee
from (A.22) we find the representations
\be
g_{ij}=
\fr1wVV'\eta_{ij}+\fr{VV'}{bq}(b_iz_j+b_jz_i)+
V^2b_ib_j+\fr1{b^2q^2}
\Bigl(VV''+(V')^2\Bigr)z_iz_j
\ee
\ses
and
\be
h_{ij}=
\fr1w VV'\eta_{ij}
+\fr1{b^2q^2}VV''z_iz_j.
\ee
The determinant of the metric tensor  is   found to read
\be
\det(g_{ij})=\fr1{w^2}
\ga
\lf(\fr1wVV'\rg)^{N-2}
V^3\det(a_{ij})
\ee
with
\be
\ga=   w^2      V''.
\ee

{

Below, the scalar $g=g(x)$ is  specified as follows:
\be
-2<g(x)<2.
\ee
We shall   apply  the convenient notation
\be
h=\sqrt{1-\fr14g^2}, \qquad
G= \fr gh.
\ee
The  {\it Finsleroid--characteristic  quadratic form}
\be
B(x,y) :=b^2+gb q+ q^2
\equiv\fr12\Bigl[(b+g_+ q)^2+(b+g_- q)^2\Bigr]>0,
\ee
where $ g_+=\frac12g+h$ and $ g_-=\frac12g-h$,
is of the negative discriminant
\be
D_{\{B\}}=-4h^2<0
\ee
and, therefore, is positively definite.

We shall use also the function $ \tau(x, w)$ defined by
 \be
 B=b^2 \tau,
 \ee
obtaining from (A.30)
 the quadratic--case representation
\be
\tau=1+g(x) w+ w^2.
\ee
We use this function to produce the  generating function $V$ according to the rule
\be
 V=\exp{\int\fr{ wd  w}{\tau }}.
\ee
Since the function (A.33) is representable in the form
\be
\tau=h^2+ \lf( w+\fr g2\rg)^2,
\ee
the integration process in (A.34) is simple, namely,
 the resultant Finslerian metric function
$
K=b V
$
(see (A.21)) is given by the following definition.

{

\ses

\ses

 {\large  Definition}.
  The scalar function $K(x,y)$ given by the formulas
\be
K(x,y)=
\sqrt{B(x,y)}\,J(x,y)
\ee
and
\be
J(x,y)=\e^{-\frac12G(x)f(x,y)},
\ee
where
\be
f=
-\arctan \fr G2+\arctan\fr{L}{hb},
\qquad  {\rm if}  \quad b\ge 0,
\ee
and
\be
f= \pi-\arctan\fr G2+\arctan\fr{L}{hb},
\qquad  {\rm if}
 \quad b\le 0,
\ee
 with
 \be
 L = q+\fr g2b,
\ee
\ses\\
is called
the {\it  Finsleroid--Finsler  metric function}.

\ses

\ses

The  function $K$ has been normalized such that
$$
0\le f\le \pi,
$$
\ses
$$
f=0,\quad {\rm if} \quad  q=0 \quad \text{ and} \quad b>0;
\qquad
f=
\pi,\quad {\rm if} \quad  q=0 \quad {\rm and} \quad b<0,
$$
and
the Finsleroid length $K(x,b^i(x))$ of the  vector $b^i$
is equal to the Riemannian length scalar $||b||=1$,
 such that
\be
K\lf(x,b^i(x)\rg)=1.
\ee

Sometimes it is convenient to use also the function
\be
A=b+\fr g2 q.
\ee
The identities
\be
L^2+h^2b^2=B, \qquad A^2+h^2 q^2=B
\ee
are valid.

The zero--vector $y=0$ is excluded from consideration.
The positive (not absolute) homogeneity holds:
$$
K(x,\la y)=\la K(x,y), \qquad \la>0, ~ \forall x, ~ \forall y.
$$

\ses

Given the function $K$ of the form (A.36), the generating function is obtained from (A.32) to read
\be
 V=\tau J
\ee
Using (A.37), it is easy to verify that
\be
\bigl(\ln V\bigr)'=\fr{ w}{\tau},
\ee
which manifests  that the integral representation (A.34) takes place.

{

\ses

\ses

 {\large  Definition}.  The arisen  space
\be
\cF\cF^{PD}_g :=\{\cR_{N};\,b_i(x);\,g(x);\,K(x,y)\}
\ee
is called the
 {\it Finsleroid--Finsler space}.

\ses

 {\large  Definition}. The space $\cR_N$ entering the above definition is called the {\it associated Riemannian space}.

\ses\ses

{\large Definition}.\, Within  any tangent space $T_xM$, the Finsleroid--metric function $K(x,y)$
 produces the {\it Finsleroid}
 \be
 \cF^{PD}_{g\,\{x\}}:=\{y\in  \cF^{PD}_{g\,\{x\}}: y\in T_xM , K(x,y)\le 1\}.
  \ee

 \ses

{\large Definition}.\, The {\it Finsleroid Indicatrix}
 $I^{PD}_{g\,\{x\}}\in T_xM$ is the boundary of the Finsleroid:
 \be
I^{PD}_{g\,\{x\}} :=\{y\in I^{PD}_{g\,\{x\}}: y\in T_xM, K(x,y)=1\}.
  \ee

 \ses

 Since at $g=0$ the  $\cF\cF^{PD}_g$--space is
Riemannian, then the body $  \cF^{PD}_{g=0\,\{x\}}$ is a unit ball and $
I^{PD}_{g=0\,\{x\}}$ is a unit sphere.

\ses

 {\large  Definition}. The scalar $g(x)$ is called
the {\it Finsleroid charge}.
The 1-form $b=b_i(x)y^i$ is called the  {\it Finsleroid--axis}  1-{\it form}.

\ses

The determinant (A.26) takes on the form
\be
\det(g_{ij})=
 V^{2N}  \fr1{ \tau^N}
\det(a_{ij}).
\ee

{

The contravariant components $g^{ij}$
of the associated Finslerian  metric tensor
can be given by the representation
\be
\fr1wVV' g^{ij}=
a^{ij}+pb^ib^j+r(b^iy^j+b^jy^i)+ty^iy^j
\ee
with
\be
r=-\fr{g}{b w},   \qquad  p+br=0,
 \qquad t=\fr g{B}(1+gw).
\ee
Therefore,
\be
\fr{K^2}{B}g^{ij}=
a^{ij}
+\fr{g}{ w}  b^ib^j
 -\fr{g}{b w} (b^iy^j+b^jy^i) +\fr g{B w}(1+g w) y^iy^j.
\ee

{

From (A.24) it follows that
\be
\fr B{K^2}  g_{ij}=
\eta_{ij}+w^2   (b_ie_j+b_je_i)+
\tau b_ib_j+
\fr{w^2}{\tau}\lf(\tau-g w\rg)e_ie_j.
\ee

In this way we obtain
    $$
\fr B{K^2}  g_{ij}=
a_{ij}
+
\fr gB\Biggl[
B w b_ib_j
-
\fr{1}{ w}v_iv_j
+
b w (b_iv_j+b_jv_i)
- b^2w^3 b_ib_j
\Biggr].
$$

Eventually, we obtain
the Finsleroid metric tensor representation
\be
\fr B{K^2}  g_{ij}=
a_{ij}
+
\fr gB\Biggl[
b^2(1+g w)  w   b_ib_j
+
b w (b_i v_j+b_j v_i)
-  \fr{1}{ w} v_i  v_j
\Biggr].
\ee

{

We can explicate the  associated vector
\be
A_k=g^{ij}A_{ijk}
\ee
 by  applying the   known general formula
$$
 A_k=K\D{\ln\lf(\sqrt{\det(g_{ij})}\rg)}
{y^k},
$$
so that
 \be
A_k=\Biggl(\ln\lf(\sqrt{\det(g_{ij})}\rg)\Biggr)'
K\D w{y^k}.
\ee

From (A.52) it follows that
\be
\fr1wVV' g^{kn}b_k= b^n -g\fr{ w}{b  \tau}y^n.
\ee

The first member of (A.22) entails the equality
\be
y_k=\fr1b K^2   \lf (b_k+\fr{b^2w^2}Be_k\rg).
\ee

With
$$
K\D w{y^k}=Vwe_k
$$
(see (A.7) and (A.19)) and (A.49),
the representation (A.56) of the vector $A_k$
is found to read
\be
A_k=   -\fr {KNg}{2B} q  e_k.
\ee
Using here (A.58), we may also write
\be
A_k= \fr{NK}2g\fr1{b w} \lf(b_k- \fr b{K^2} y_k\rg).
\ee
Another convenient form is
\be
A_k=\fr {NK}2g\fr1{ qB}(q^2b_k- bv_k).
\ee

{

From (A.57) and (A.60) it follows immediately that
\be
A^k=\fr{Ng}{2Kb w}
\Bigl[Bb^k
- b(1+g w)y^k
\Bigr].
\ee
\ses
It can readily be seen that the representations (A.60) and (A.62) entail
\be
A^kA_k=\fr{N^2g^2}4.
\ee
\ses
We have
\be
y_kb^k= b\fr{K^2}B(1+g w),
\ee
\ses
\be
A_kb^k  =
\fr {KNg}{2B} b w, \qquad
b_kA^k=\fr{Ng}{2K}b w,
\ee
and
\be
Kg^{kj}b_j=
 \fr{2b w}{Ng}A^k
 +  b  l^k,
\ee
\ses
together with
\be
K^2g^{kj}b_j=
  Bb^k
- g b  w   y^k.
\ee

The relation (A.59) can be inverted, yielding
\be
e_k=  -\fr {2B}{KNgq} A_k.
\ee
Taking into account the formulas (A.7), (A.19),  and (A.68), we may write simply
\be
\D{ w}{y^i}= - \fr{2B}{b^2} \fr1 {KNg } A_i.
\ee
This formula is convenient to use in many involved evaluations.

     \ses

With (A.52), it follows that
\be
Kg^{kj}g_j=
\fr1KB\lf(g^k-(bg)b^k\rg)
+  \fr{2b w}{Ng}(bg)   A^k
+   \fr{2}{N}[b(bg)-(yg)]   A^k+b(bg)l^k.
\ee

{

Now we perform differentiation of the functions indicated in (A.29)
with respect to the Finsleroid parameter $g$, obtaining
 \be
  \D hg= -\fr14 G, \qquad  \D Gg= \fr1{h^3}, \qquad \D{\lf(\fr Gh\rg)} g =\fr1{h^4} \lf(1+\fr{g^2}4 \rg),
  \ee
\ses
\be
\D fg= -\fr1{2h} +  \fr{b}B\Bigl(\fr14G  q +   \fr1{2h}b\Bigr),
\ee
and
 \be
2K  \D Kg=b  qJ^2 - \fr1{h^3}fK^2+
G \Biggl[\fr1{2h}-\fr{b}B\Bigl(\fr14G  q +   \fr1{2h}b\Bigr)\Biggr]K^2,
  \ee
\ses
or
\be
 \D {K^2}g=M
K^2,
\ee
where
$$
M=
 \fr{b q}B - \fr1{h^3}f+
\fr12\fr{G}{hB} ( q^2+\frac12 gb q),
$$
\ses
or
\be
M=
 - \fr1{h^3}f+
\fr12\fr{G}{hB}   q^2+  \frac1{h^2B} b q.
\ee

{

Using
$$
e_i=-b_i+\fr b{q^2}v_i,
$$
we find
\be
M_i=\fr{4  q^2}{gNBK} A_i.
\ee

\ses

Differentiating (A.74) with respect to $y^i$ just yields
\be
 \D {{y_i}}g=
My_i
+\fr 12 M_iK^2.
\ee
In view of (A.59) and  (A.76), we may write
\be
 \D {{y_i}}g=
My_i
-\fr{q^2 qK^2}{B^2}e_i,
\ee
or
\be
 \D {{y_i}}g=
My_i
+\fr{2  q^2K}{gNB} A_i.
\ee

{

\nin
With
the tensor
\be
{\cal H}_{ij}=h_{ij}-\fr{A_iA_j}{A_nA^n},
\ee
 we can arrive at
\be
 \D {{g_{ij}}}g=
Mg_{ij}
+\fr{ q^2}{B}   \fr 2{gN}  (A_il_j +l_iA_j )
- \fr {b  q}B {\cal H}_{ij}
-\fr{2b  q}{B}  \fr 4{g^2N^2} A_i   A_j.
\ee

From (A.61) we obtain the simple result
\be
\D{A_i}g=
 \Biggl(  \fr12M        +\fr1g   -  \fr bB  q       \Biggr)
A_i.
\ee

From (A.62) it follows that
$$
\D{\lf(   \fr K{Ng}     A^k \rg)} g=  -  \fr{1}{2}   \Bigl(y^k-   b b^k     \Bigr).
$$

From (A.75) we get
$$
\D Mg= \fr{3 g}{4h^2} M
+\fr1{h^2} \fr{  q^2 }{B}-
   \fr1{h^2} \fr{  q^2 }{B^2}
  \lf(b^2+\fr12g b  q \rg).
$$

{

Next,  with (A.62)
we now evaluate the derivative
\be
\D{A^k}{y^n}  =  -      \fr {1}{K} A_n l^k
-\fr{Ng}{2K w}   (1+g w)\lf( \de^k{}_n-l^kl_n-\fr4{N^2g^2}A^kA_n\rg)
-\fr2N \fr1K A^kA_n.
\ee

\ses

The previous representation can be
 written as
\be
\D{A^k}{y^n}  =  -      \fr {1}{K} A_n l^k
-\fr{Ng}{2K w}   (1+g w){\cal H}^k{}_n
-\fr2N \fr1K A^kA_n
\ee
in terms of the tensor
\be
{\cal H}^k{}_n= \de^k{}_n-l^kl_n-\fr4{N^2g^2}A^kA_n
\ee
given by (A.80).

The Cartan tensor
\be
A_{ijk} =\fr12K\D{g_{ij}}{y^k}
\ee
 takes on the form
\be
A_{ijk}=  \fr1N\lf(h_{jk}A_i+h_{ik}A_j+h_{ij}A_k   -  \fr1{A^hA_h}A_iA_jA_k  \rg),
\ee
or in terms of the tensor (A.85),
\be
A_{ijk}=  \fr1N\lf({\cal H}_{jk}A_i+ {\cal H}_{ik}A_j+ {\cal H}_{ij}A_k  +  \fr2{A^hA_h}A_iA_jA_k  \rg).
\ee
We can conclude that
\be
A^kA_{ijk}=  \fr1N \lf(A_iA_j+h_{ij}A_kA^k\rg) \equiv  \fr1N \lf(2A_iA_j+{\cal H}_{ij}A_kA^k\rg).
\ee
Also, we find
\be
A_{ijk}A^k{}_{mn}= \fr1N A_i  A_{jmn}  + \fr1N A_j  A_{imn} +   \fr2{N^2}{\cal H}_{ij}  A_mA_n
+\fr1N \fr1N A_kA^k  {\cal H}_{ij}   {\cal H}_{mn}.
\ee

\ses

{

From (A.87) it follows that
\be
A^{ijk}
A_{ijk}=
\fr1N \Bigl(3-\fr2N\Bigr)
A^hA_h.
\ee

\ses

\ses

The curvature of  indicatrix is well--known to be described by the tensor
\be
\hat R_i{}^j{}_{mn} := \fr1{K^2}\Bigl(\3Ahjm\3Aihn-\3Ahjn\3Aihm\Bigr).
\ee
Inserting (A.90) in (A.92), we find that
\be
K^2\hat R_{ijmn}
=
\fr1{N^2}  (A^kA_k)
\Bigl(
 h_{in}h_{mj}  -h_{im}h_{nj}
\Bigr).
\ee
\ses
Contracted objects
\be
R_{im}=\hat R_i{}^j{}_{mj}
\ee
and
\be
R=\hat R_i{}^j{}_{mj}g^{im}
\ee
are found to be
\be
R_{im}
=
\fr1{K^2}  \lf(\fr2N-1\rg)    \fr1N(A^kA_k)
h_{im}
\ee
and
\be
R
=
\fr1{K^2}  \lf(\fr2N-1\rg)  (N-1)   \fr1N (A^kA_k).
\ee

{

If we consider the derivative
$$
\D{A_k}{y^n}=  g_{km}\D{A^m}{y^n} +\fr2K A_{kmn}A^m
$$
and apply (A.84) together with (A.89), we obtain
\be
\D{A_k}{y^n}  =  -      \fr {1}{K} A_n l_k
-\fr{Ng}{2K w}  {\cal H}_{kn}
+\fr2N \fr1K A_kA_n.
\ee

Also,
\be
\tau_{ij}~:
=
A_{ij}-A_k\3Aikj,
\ee
where
\be
A_{ij}~:= K \D{ A_i}{ y^j}   + l_iA_j.
\ee
The identities
$$
\tau_{ij}y^j=0,      \qquad        \tau_{ij}A^j=0
$$
hold. We obtain
\be
\tau_{ij}=-\fr N4\fr{g(2b+g q)}{ q}{\cal H}_{ij}.
\ee
\ses

Differentiating (A.85) leads  to
\be
K\D{{\cal H}^k{}_n}{y^m}= -{\cal H}^k{}_ml_n -l^k{\cal H}_{nm}
+\fr2{Ng} \fr1{ w}  {\cal H}_{nm}A^k
+\fr2{Ng}   \fr1{ w}  (1+g w){\cal H}^k{}_m  A_n.
\ee

{

If we apply (A.98) to (A.87), we get
$$
\D{A_{ijk}}{y^n}
= \fr1K    \fr2N\lf(A_{jkn}A_i+A_{ikn}A_j+A_{ijn}A_k  \rg)
- \fr1K \Bigl(l_jA_{kni}+l_iA_{knj}+l_kA_{ijn}\Bigr)
$$

\ses

$$
+  \fr1K \fr1N \fr2N \Bigl[{\cal H}_{jk} A_iA_n  +{\cal H}_{ik} A_jA_n+{\cal H}_{ij} A_kA_n\Bigr]
$$

\ses

\be
 -  \fr{g}{2K w}
  \Bigl[ {\cal H}_{jk}{\cal H}_{in} + {\cal H}_{ik}{\cal H}_{jn}+ {\cal H}_{ji}{\cal H}_{kn}\Bigr].
\ee
\ses
If we introduce the tensor
\be
\tau_{ijkn}~:=K\D{A_{ijk}}{y^n}
+  l_jA_{kni}+l_iA_{knj}+l_kA_{ijn}
- A_{jkm}A^m{}_{in}  - A_{ikm}A^m{}_{jn}  - A_{jim}A^m{}_{kn},
\ee
then we can conclude from (A.103) and
 (A.90) that
\be
\tau_{ijkn}=
 -  \fr{g(2+g w)}{4 w}
  \Bigl[ {\cal H}_{jk}{\cal H}_{in} + {\cal H}_{ik}{\cal H}_{jn}+ {\cal H}_{ji}{\cal H}_{kn}\Bigr].
\ee

{

In calculating the Finsleroid--case spray coefficients
there appears the  vector
\be
E^k~:= M (yg)y^k
+\fr12  K^2 (yg) M_h g^{kh}
-\fr12 MK^2g_hg^{kh}.
\ee
Noting (A.76), we obtain
\be
E^k = M (yg)y^k
+  K \fr{2 b^2 w^2}{gNB}(yg)  A^k
-\fr12 MK^2g_hg^{kh}.
\ee

We may calculate the derivative
\be
E^k{}_n~:=\D{E^k}{y^n}.
\ee
We obtain
$$
E^k{}_n = MT^k{}_n
+  K \fr{2 b^2 w^2}{gNB} A^kg_n
$$

\ses

$$
+   \fr{2 b^2 w^2}{gNB}(yg)  A^k  l_n
-   \fr{4  w}{gN}(yg)   \fr {2}{Ng} A_n   A^k
$$

\ses

\ses

\be
+   \fr{4 b^2 w^3}{gNB}(yg)   \fr {2}{Ng}  A_n
 A^k
+   \fr{2 b^2 w^2}{gNB}(yg)
\Biggl[
       A_n l^k
-\fr{Ng}{2 w}   (1+g w){\cal H}^k{}_n
\Biggr]
-\fr12 K^2M_n      g_hg^{kh}
\ee
\ses
with
$$
  T^k{}_n=g_ny^k+  (yg)\de^k{}_n   -  \fr12  Kg_hg^{kh}l_n
  -\fr12K   \D{\Bigl(Kg_hg^{kh}\Bigr)}{y^n}.
  $$
The contraction $Kg_hg^{kh}$ is to be taken from  (A.70).

{

Using (A.81) yields
$$
\lf(\fr{\partial g_{kj}}{\partial g} g_i
+\fr{\partial g_{ik}}{\partial g}  g_j
-\fr{\partial g_{ij}}{\partial g}  g_k
\rg)A^k
=
$$

\ses

      \ses

\ses

$$
\Biggl(MA_j   +\fr{ q^2}{B}   \fr {gN} 2  l_j   -\fr{2b q}{B}     A_j   \Biggr)g_i
+ \Biggl(MA_i   +\fr{ q^2}{B}   \fr {gN} 2  l_i   -\fr{2b q}{B}     A_i   \Biggr)g_j
$$

\ses

\ses

\be
-\Biggl(Mg_{ij}   +\fr{ q^2}{B}   \fr 2{gN}  (A_il_j +l_iA_j )   - \fr {b q}B {\cal H}_{ij}
-\fr{2b q}{B}  \fr 4{g^2N^2} A_i   A_j
\Biggr)
g_kA^k.
\ee

\ses

{

We can write down  simple explicit representations for the partial
 derivatives with respect to $x$,
using the notation
\be
 g_j=\D g{x^j},  \qquad s_k   =y^m\nabla_kb_m.
 \ee
We may use the equality
\be
\D {w}{x^k}=
-\fr 1{b^2q}S^2s_k+\De
\ee
(see (A.6)),
 where   $\De$ symbolizes the summary of the terms which involve partial derivatives
of the input Riemannian metric tensor $a_{ij}$ with respect to the coordinate variables $x^k$.
We use
 the Riemannian covariant derivative
\be
\nabla_ib_j~:=\partial_ib_j-b_ka^k{}_{ij},
\ee
where
\be
a^k{}_{ij}~:=\fr12a^{kn}(\prtl_ja_{ni}+\prtl_ia_{nj}-\prtl_na_{ji})
\ee
are the
Christoffel symbols given rise to by the associated Riemannian metric.

First of all, we differentiate the quadratic form $B$ given by (A.30), obtaining
\be
\D{B}{x^j}=b  q  g_j +\fr2b(B-S^2)  s_j
-  g\fr1{ q}     S^2   s_j
+\De.
\ee

Also,
starting with (A.21),
we find
\be
\D {K}{x^j}= \fr12 M Kg_j +  K\fr b{B} g w s_j
+\De,
\ee
where we have used (A.74).

Eventually, the following sufficiently simple representation is obtained:
 $$
 \D{A_i}{x^j}=   \Biggl(  \fr12M        +\fr1g   -  \fr bB  q       \Biggr)A_i g_j
+
\fr {NK}2g\fr1{ qB}S^2 \nabla_jb_i
$$

\ses

\ses

$$
+ \fr b{B} g w s_j   A_i
 -\fr1B  \Biggl(\fr2b(B-S^2)  s_j
-  g\fr1{ q}     S^2   s_j
\Biggr)
 A_i
 $$

\ses

\ses

\be
+\fr1{ w^2}
  \fr 1{b^3}S^2s_j A_i
-\fr {NK}2g\fr1{ qB}  \fr 1{b}S^2
s_j  b_i +\De.
\ee

{

        \bigskip

\setcounter{equation}{0}

\bc
{ \bf Appendix B: Finsleroid--Finsler   spray coefficients }
\ec

\ses

\ses

Evaluations involve  the induced {\it spray coefficients}
\be
 G^k=\ga^k{}_{ij}y^iy^j,
\ee
which entail the  coefficients
\be
 G^i{}_k:~=\D{ G^i}{y^k}, \qquad  G^i{}_{km}:~=\D{ G^i{}_k}{y^m},
\qquad  G^i{}_{kmn}:~=\D{ G^i{}_{km}}{y^n},
\ee
and
\be
\bar G^i=\fr12G^i, \qquad \bar  G^i{}_k =\fr12 G^i{}_k, \qquad  \bar G^i{}_{km}  =\fr12G^i{}_{km},
\qquad  \bar G^i{}_{kmn}=\fr12G^i{}_{kmn}.
\ee
The homogeneity gives rise to the identities
\be
2G^i= G^i{}_ky^k, \qquad  G^i{}_k=G^i{}_{km}y^m,
\qquad  G^i{}_{kmn}y^n=0.
\ee

The pair $(x,y)$, --- the so--called {\it line element}, --- is the argument of the Finslerian objects;
\be
\ga^k{}_{ij}~:=\fr12 g^{kn}\lf( \D{g_{ni}}{x^j}  +  \D{g_{nj}}{x^i}
-\D{g_{ji}}{x^n}\rg)
\ee
are the
Christoffel symbols given rise to by the Finsleroid--Finsler   metric function $K$.

{

Below, the abbreviation $h$ means {\it horizontal}.
On the basis  of the  above  coefficients
the Finslerian
{\it connection coefficients}
$\Ga^k{}_{ij}$  of $h$--type
are constructed according to the well-known  rule:
\be
\Ga^k{}_{ij}=\ga^k{}_{ ij}-\bar G^n{}_i\3Cnkj-\bar G^n{}_j\3Cnki+\bar G^{kn}C_{nij}
\ee
with
\be
\bar G^n{}_i=\ga^n{}_{ ij}y^j-2\bar G^m\3Cmni=\Ga^n{}_{ij}y^j=\fr12 G^n{}_i
\ee
and
\be
2\bar G^m=\ga^m{}_{ij}y^iy^j=
\bar G^m{}_iy^i=\Ga^m{}_{ ij}y^iy^j=G^m,
\ee
where
 $C_{nij}=A_{nij}K^{-1}$  and $\3Cnkj=\3Ankj K^{-1}$.
By the help of these coefficients
 the $h$-{\it covariant derivatives} of tensors are constructed
as exemplified by
\be
A_{i|j}~:= \prtl_jA_i-\bar G^k{}_j\D{A_i}{y^k}-\Ga^k{}_{ij}A_k
\ee
(see [1,2]).

For   the $hh$-curvature tensor $R^i{}_k$ we use the formula
\be
K^2R^i{}_k~:=
2\D{\bar G^i}{x^k}-\D{\bar G^i}{y^j}\D{\bar G^j}{y^k}
-y^j\Dd{\bar G^i}{x^j}{y^k}
+2\bar G^j\Dd{\bar G^i}{y^k}{y^j}
\ee
(which is tantamount to the definition (3.8.7) on p. 66 of the book [2]).
The concomitant tensors
\be
R^i{}_{km}~:=\fr1{3K}\Biggl(\D{(K^2R^i{}_{k})}{y^m}-\D{(K^2R^i{}_{m})}{y^k}
\Biggr),
\ee
\ses
and
\be
R_n{}^i{}_{km}~:=\D{(KR^i{}_{km})}{y^n}
-\Bigl( \dot A^i{}_{nm|k} - \dot A^i{}_{nk|m}
+\dot A^i{}_{uk}\dot A^u{}_{nm} -  \dot A^i{}_{um}\dot A^u{}_{nk} \Bigr)
\ee
(see p. 60 in the book [2])
arise.
We have
\be
R^i{}_ky^k=0
\ee
and
\be
R^i{}_{km}y^m=KR^i{}_k.
\ee

{

In calculations, it proves convenient to write the derivative (B.9) in the alternative form

\be
A_{i|j}=\prtl_jA_i-\bar G^k{}_j\tau_{ik} \fr1K -    {\wt\Ga}^k{}_{ij}  A_k+\bar G^k{}_jA_kl_i \fr1K
\ee
($l_i=g_{ik}l^k$)
with
\be
{\wt\Ga}^k{}_{ij}=\ga^k{}_{ ij}-\bar G^n{}_i\3Ankj  \fr1K +\bar G^{kn}A_{nij}  \fr1K
\ee
and
\be
\tau_{ij}
=
A_{ij}-A_k\3Aikj,
\ee
where
\be
A_{ij} = K\D{ A_i}{ y^j}   +l_iA_j.
\ee
We  know that
\be
\tau_{ij}=-\fr N4\fr{g(2b+g  q)}{ q}{\cal H}_{ij}
\ee
(see (A.101)).

\ses

{

By means of attentive (lengthy) evaluations we can arrive at   the following assertion.

\ses

\ses

THEOREM B1.  {\it The explicit form of the
 spray coefficients
of the  Finsleroid--Finsler space reads}
\be
G^k=
  g q
\Bigl[a^{kj}+(pb^k+ry^k)b^j\Bigr]
y^h(\nabla_hb_j- \nabla_jb_h)
+\fr{g}{q}
\lf(y^k-bb^k\rg)(ys)
+a^k{}_{mn}y^my^n
+E^k,
\ee
\ses
where
$p$ and $r$ are the quantities presented in (A.51)
and
$E^k$ is the vector (A.107).
The notation
\be
(ys)=y^hs_h
\ee
has been used, where
\be
s_k      =y^m\nabla_kb_m.
\ee


A careful consideration of the formulas (B.20)--(B.22) shows that the following theorem is valid.

\ses\ses

 THEOREM B2.
{\it The Finsleroid--Finsler space ${\mathbf\cF^{PD}_g } $
 is of the Landsberg type if and only if
the following three conditions hold:
the Finsleroid charge is a constant
\be
g=const,
\ee
the input $1$-form $b$ is closed
\be
\partial_ib_j-\partial_jb_i=0,
\ee
and the expansion
\be
\nabla_mb_n=
k\lf(a_{mn}-b_mb_n\rg)
\ee
takes place, where
$
k=k(x)
$
 is a scalar.
}

\ses\ses

Under the conditions of this theorem,
we have $E^k=0$ and $(ys)=kq^2$, so that  the representation (B.20) reduces to
\be
G^k=
gkq \lf(y^k-bb^k\rg)
+a^k{}_{mn}y^my^n.
\ee

\ses

It can readily be seen that
at any dimension $N\ge3$ the Berwald case corresponds to $ k=0$,
so that
\be
 G^m_{\text{\{Berwaldian\}}}=    a^m{}_{jn}y^jy^n.
\ee
At the dimension $N=2$,
 the Finsleroid--Finsler space ${\mathbf\cF^{PD}_g } $
 is of the Landsberg type if and only the space is of the Berwald type
(independently of  value of $k$).

\ses

\ses

NOTE.   Theorem B2 is known from the previous publications [4-7].
Theorem B1 is a new result. When $g=const$, the coefficients $E^k$   vanish identically,
in which case the above spray coefficients (B.20) coincide with the spray coefficients given by Eq. (4.5)
in [7].

{

            \bigskip

\setcounter{equation}{0}

\bc
{\bf  Appendix C:  Finsleroid--involutive  tensors}
\ec

\ses\ses

\ses

Henceforth, we assume the  involutive case
\be
g_i=\mu b_i, \qquad \mu =\mu(x),
\ee
which entails
\be
b(bg)=(yg),   \qquad (bg)=\mu, \qquad (yg)=\mu b.
\ee
The notation $(yg)=y^ig_i$ and $(bg)=b^ig_i$ is used; $g_i=\partial g/\partial x^i$.
Under these conditions, the formula (A.70) reads merely
\be
Kg^{kj}g_j=
   \fr{2b  w}{Ng}(bg)   A^k
 +b(bg)l^k
\ee
and the representation (A.107) reduces to read
\be
E^k =  \fr12 M (yg)y^k
-  \wh   MK  \fr1{Ng}    w  (yg)   A^k,
\ee
\ses
where
\be
\wh M = M  - \fr{2 b^2 w}{B}.
\ee
Differentiating this scalar yields the simple result:
\be
\D{\wh M}{y^m}= \fr{4 b^2}{B}       \fr1 {KNg } A_m.
\ee
\ses

The eventual representation reads
\be
E^k{}_n = MT^k{}_n
+   \fr{4 b^2 w^2}{gNB}(yg)  A^k  l_n
-   \fr{4  w b^2(1+g w)}{gNB}  \fr {2}{Ng}(yg)  A_n   A^k
-   \fr{ b^2 w}{B}(yg)
(1+g w){\cal H}^k{}_n
\ee
\ses
with
$$
  T^k{}_n= (yg)l^kl_n +   \fr2{Ng}  w   (yg)  l^k   A_n
  +   \fr12 (yg){\cal H}^k{}_n
-    \fr 2{Ng}    w  (yg)   A^k
   l_n
$$

 \ses

\ses

\be
+    \fr 2{Ng}    \fr {2(1+g w)}{Ng}    (yg)  A_n  A^k
+  \fr12(yg)   (1+g w){\cal H}^k{}_n.
\ee

{

We find
\be
A_kE^k = - \fr {Ng}4 K     w (yg)
\Bigl(M  - \fr{2 b^2 w}{B}  \Bigr)
\ee
\ses
and
$$
A_kE^k{}_n = M  \lf(   -    \fr {Ng}2      w
   l_n
+    (1+g w)     A_n
\rg)
(yg)
+   \fr{2 b^2 w^2}{B}  \fr {Ng}2 (yg) l_n
-   \fr{2  w b^2(1+g w)}{B}     (yg)  A_n,
$$
or
\be
A_kE^k{}_n = \Bigl(M  - \fr{2 b^2 w}{B}  \Bigr)
 \lf(   -    \fr {Ng}2      w
   l_n
+    (1+g w)     A_n
\rg)
(yg).
\ee
Also,
\be
E^k{}_nA^n = \Bigl(M  - \fr{2 b^2 w}{B}  \Bigr)
    (1+g w)  (yg)   A^k
+M    \fr {Ng}2      w
(yg)   l^k
\ee
\ses
and
\be
{\cal H}_{km}E^k{}_n = M
\Biggl[   \fr12 (yg) {\cal H}_{mn}  +  \fr12(yg)   (1+g w)  {\cal H}_{mn}  \Biggr]
-   \fr{ b^2 w}{B}(yg)
(1+g w)   {\cal H}_{mn},
\ee
together with
\be
E^k{}_n {\cal H}^n{}_i= M
\Biggl[   \fr12 (yg) {\cal H}^k{}_i  +  \fr12(yg)   (1+g w)  {\cal H}^k{}_i  \Biggr]
-   \fr{ b^2 w}{B}(yg)
(1+g w)   {\cal H}^k{}_i.
\ee

\ses

{

We obtain
$$
-E^n{}_i\3Ankj  + E^{kn}A_{nij} =
-
\fr1N
\Biggl[
M \fr12 (yg) {\cal H}_{ij}
-  \wh M
 \fr12(yg)   (1+g w)  {\cal H}_{ij}
\Biggr]
A^k
$$

\ses

$$
+
\wh M
   \fr {g}2      w
(yg) {\cal H}^k{}_j l_i
+
\wh M
 g      w
(yg)   \fr 1{A^hA_h}A^kA_j l_i
$$

\ses\ses

$$
+
\fr1N
\Biggl[
M \fr12 (yg) {\cal H}^k{}_j
-\wh M
 \fr12(yg)   (1+g w)  {\cal H}^k{}_j
\Biggr]
A_i
$$

\ses

\be
+   M    \fr {g}2     w
(yg) {\cal H}_{ij}   l^k
+M  g      w
(yg)   \fr 1{A^hA_h}A_iA_j l^k.
\ee

\ses

The last formula just entails
\be
\Bigl(-E^n{}_i\3Ankj  + E^{kn}A_{nij} \Bigr)A_k=
-
\fr{Ng^2}8
\Bigl[
M (yg) {\cal H}_{ij}
-\wh M
 (yg)   (1+g w)  {\cal H}_{ij}
    \Bigr]
+
\wh M
 g      w
(yg)  A_j l_i.
\ee

The formula (A.110) can be written as
$$
 \lf(\fr{\partial g_{kj}}{\partial g} g_i
+\fr{\partial g_{ik}}{\partial g}  g_j
-\fr{\partial g_{ij}}{\partial g}  g_k
\rg)A^k
=
$$

\ses

\ses

\be
\mu  \wh M   \Biggl(
\fr bK(A_jl_i+A_il_j)     -\fr{Ng}{2K}  q  l_il_j
 -\fr{Ng}{2K}  q  {\cal H}_{ij}
+\fr2{NKg}b wA_iA_j
 \Biggr).
\ee

      \ses

\ses

{

Let us assume that the $b$-parallel case
$$
\nabla_ib_j=0
$$
takes place.
Then
 it proves convenient to write the derivative (B.15) in the form
\be
A_{i|j}=\prtl_jA_i
-\fr12 E^k{}_j\tau_{ik} \fr1K - {\wt\Ga}^{ k}{}_{ij} A_k+ \fr12  E^k{}_jA_kl_i \fr1K
+\De
\ee
with
\be
{\wt\Ga}^{ k}  {}_{ij}=\ga^k{}_{ ij}-\fr12 E^n{}_i\3Ankj  \fr1K +\fr12 E^{kn}A_{nij}  \fr1K
+\De,
\ee
and   required insertions lead  to the following simple result:
\ses\\
\be
K A_{i|j}=    \fr1g KA_ig_j
+  \fr14  \mu  \wh M     \fr {Ng}2    \fr{B} {  q}     {\cal H}_{ij}
+\fr 14   \mu M  \fr{Ng}{2}
\Biggl(    \fr{b(2b+g  q)}{2  q}   +   q +gb   \Biggr)       {\cal H}_{ij}.
\ee

{

Next, we consider the derivative tensor
\be
E^k{}_{nm}~:=\D{E^k{}_n}{y^m}
\ee
to find the contraction
$$
A^nE^k{}_{nm} =\D{A^nE^k{}_n}{y^m}- E^k{}_n\D{A^n}{y^m}.
$$

Make required cancellation and
use (C.4),
obtaining
$$
A^nE^k{}_{nm} =
 \wh M [ w-g(1+g w)]
 \fr2{NKg}  (yg)  A_m A^k
$$

\ses

$$
+\Biggl[\wh M (1+g w) A^k  +M \fr {Ng}2   w  l^k
\Biggr]
\fr 1{K} (yg)l_m
$$

\ses

$$
+
      2  \fr2 {KNg }
  (yg)   A^k  A_m
-\wh M
    (1+g w)  (yg)
      \fr {1}{K} A_m l^k
$$

\ses

$$
+ \fr{2  q^2}{BK}     w
(yg)   l^k A_m
+M    \fr {Ng}2  (yg)    w \fr1K {\cal H}^k{}_m
+  M   \fr{Ng}{4K w}   (1+g w)
  (yg) {\cal H}^k{}_m.
$$
We get
$$
E^nE^k{}_{nm} =\fr12 M (yg) E^k{}_m
-
\wh M
K  \fr1{Ng}    w  (yg)   A^n E^k{}_{nm}.
$$

Now, write (C.7)--(C.8) as follows:
$$
E^k{}_n = M \Biggl((yg)l^kl_n +   \fr2{Ng}   w   (yg)  l^k   A_n
  +   \fr12 (yg){\cal H}^k{}_n
  \Biggr)
  $$

\ses

\be
+\wh M\Biggl[-    \fr 2{Ng}     w  (yg)   A^k     l_n
+    \fr 2{Ng}    \fr {2(1+g w)}{Ng}    (yg)  A_n  A^k
+  \fr12(yg)   (1+g w){\cal H}^k{}_n
\Biggr],
\ee

{

\nin
or
$$
-E^k{}_nE^n{}_m
+
2E^nE^k{}_{nm}
=    \fr14 M^2(yg)^2      {\cal H}^k{}_m
$$

\ses

$$
+  M\wh M(yg)^2
\Biggl(
 \fr2{Ng}\fr2{Ng}  \fr B{b^2}      A^kA_m
-  \fr14   (1+gw){\cal H}^k{}_m
\Biggr)
$$

\ses

$$
-\wh M^2   (yg)^2
\Biggl(
     \fr {4}{N^2g^2}  (1+gw)^2   A^kA_m
+  \fr14   \fr B{b^2}    {\cal H}^k{}_m
\Biggr)
$$

\ses

\ses

\be
-
\wh M
  \fr 2{Ng}    w  (yg)^2
\Biggl[
       2\fr2 {Ng }       A^k  A_m
-\wh M
    (1+gw)
l^k A_m
+ \fr{2 \wt q^2}{B}    w
   l^k A_m
+M    \fr {Ng}2         w  {\cal H}^k{}_m
\Biggr].
\ee

{


Below we again assume
 that the $b$-parallel condition
$$
\nabla_ib_j=0
$$
holds.
Since
$$
\D{M}{x^m}=\D M g  g_m
+\De,
$$
we can straightforwardly come   from (C.4)    to
$$
\D{E^k}{x^m}=  \fr1{\mu} E^k \mu_m
+  \fr12 (yg)y^k \D Mg g_m
 - \fr K{Ng}    w     A^k (yg) \D{  \wh   M}g g_m
 + \fr12 \wh   M    w  (yg)    (y^k-   b b^k ) g_m
 +\De,
$$
where $\mu_m=\partial \mu/\partial x^m$
and the formula placed below (A.82) has been applied.
In this way we obtain
\ses
$$
y^m\D{E^k}{x^m}     = \fr1{\mu}(y\mu)E^k
+  \fr12 (yg)^2 y^k \D Mg
 -  \fr K{Ng}    w     A^k (yg)^2 \D{  \wh   M}g
 + \fr12 \wh   M    w  (yg)^2  ( y^k-   b b^k )
 +\De
$$
\ses
and
\be
2\D{E^i}{x^k} -y^j\Dd{E^i}{x^j}{y^k}
=  2\fr1{\mu} E^i \mu_k    - \fr1{\mu}(y\mu)E^i{}_k    +(yg)^2S^i{}_k
+\De.
\ee

\ses

Now it is easy to continue the calculation:
we shall use the equality
\be
 \D{  \wh   M}g= \D{     M}g +\fr{2b^2 q^2}{B^2}
 \ee
ensued from (C.5).

We find
$$
S^i{}_k=
-\fr12  \Biggl( \fr{B}{b^2}  \fr 2{Ng}\fr2 {Ng } A_k     A^i
 +  (1+gw){\cal H}^i{}_k
+  {\cal H}^i{}_k     +   (1+gw +w^2)  \fr 2{Ng} \fr2{Ng} A^iA_k
\Biggr)
 \D{     M}g
$$

\ses
\ses

$$
-\Biggl( \fr{B}{b^2} \fr 2{Ng}\fr2 {Ng } A_k     A^i
+   (1+gw){\cal H}^i{}_k     +   (1+gw +w^2) \fr 2{Ng}  \fr2{Ng} A^iA_k
\Biggr)
\fr{b^2q^2}{B^2}
$$

\ses
\ses

\be
- \wh   M    w  \Biggl(
\fr12 {\cal H}^i{}_k
+  \fr 2{Ng }  \fr 2{Ng }  A^i A_k
 - \fr2{Ng}  w   A_kl^i
\Biggr),
\ee
where we  must  insert the derivative
\be
\D Mg= \fr1{h^2}
\Biggl[
\fr{3 g}{4} M
+ \fr{q^2 }{B}-
  \fr{q^2 }{B^2}
  \lf(b^2+\fr12g bq \rg)
\Biggr]
\ee
(see the formulas below (A.82)).

{

The respective {\it involutive curvature tensor} $R^{i}{}_k$ is constructed according to
\be
K^2 R^{i}{}_k=
2\D{\bar E^i}{x^k} -y^j\D{\bar E^i{}_k}{x^j}
-\bar E^i{}_n\bar E^n{}_k
+
2\bar E^n\bar E^i{}_{nk}
+y^na_n{}^i{}_{km}y^m,
\ee
where
\be
\bar E^i=\fr12E^i, \qquad \bar  E^i{}_k =\fr12 E^i{}_k, \qquad  \bar E^i{}_{nk}  =\fr12E^i{}_{nk},
\ee
and $a_n{}^i{}_{km}$ stands for the Riemannian curvature tensor of the associated Riemannian space.
The explicit formulas (C.22)--(C.26) must be inserted in  (C.27),  yielding
the following result:
$$
K^2 R^{i}{}_k=
 \fr1{\mu} E^i \mu_k    - \fr1{2\mu}(y\mu)E^i{}_k
$$

\ses
\ses

$$
   +\fr12   (yg)^2
\Biggl[
-\fr12  \Biggl(2 \fr{B}{b^2}  \fr 2{Ng}\fr2 {Ng } A_k     A^i
 +  (1+gw){\cal H}^i{}_k
\Biggr)
 \D{     M}g
$$

\ses
\ses

$$
-\Biggl(2 \fr{B}{b^2} \fr 2{Ng}\fr2 {Ng } A_k     A^i
+   (1+gw){\cal H}^i{}_k
\Biggr)
\fr{b^2q^2}{B^2}
$$

\ses
\ses

$$
- \wh   M    w  \Biggl(
\fr12 {\cal H}^i{}_k
+  \fr 2{Ng }  \fr 2{Ng }  A^i A_k
 - \fr2{Ng}  w   A_kl^i
\Biggr)
\Biggr]
$$

\ses

\ses

              \ses

\ses

$$
+    \fr1{16}     M^2(yg)^2      {\cal H}^k{}_m
+ \fr14 M\wh M(yg)^2
\Biggl(
 \fr2{Ng}\fr2{Ng}  \fr B{b^2}      A^kA_m
-  \fr14   (1+gw){\cal H}^k{}_m
\Biggr)
$$

\ses

\ses

$$
-  \fr14  \wh M^2   (yg)^2
\Biggl(
     \fr {4}{N^2g^2}  (1+gw)^2   A^kA_m
+  \fr14   \fr B{b^2}    {\cal H}^k{}_m
\Biggr)
$$

\ses

\ses

$$
-
 \fr14      \wh M
  \fr 2{Ng}    w  (yg)^2
\Biggl(
       2\fr2 {Ng }       A^k  A_m
-\wh M
    (1+gw)
l^k A_m
+ \fr{2 q^2}{B}    w
   l^k A_m
+M    \fr {Ng}2         w  {\cal H}^k{}_m
\Biggr)
$$

\ses

\ses

\be
+y^na_n{}^i{}_{km}y^m.
\ee

{

In evaluations, we apply the representation
\be
\D{M_i}g= -  4 \fr{bq^3 }{B^2}
\fr {2}{KNg}  A_i.
\ee

{

          \bigskip

\ses\ses

\def\bibit[#1]#2\par{\rm\noindent\parskip1pt
                     \parbox[t]{.05\textwidth}{\mbox{}\hfill[#1]}\hfill
                     \parbox[t]{.925\textwidth}{\baselineskip11pt#2}\par}

\bc {\bf  References} \ec
\bigskip

\bibit[1] H. RUND,  The Differential Geometry of Finsler spaces,
{\it Springer, Berlin,} 1959.

\bibit[2]  D. BAO, S. S. CHERN, and Z. SHEN:   An
Introduction to Riemann-Finsler Geometry,
{\it  Springer, New York, Berlin,}  2000.

\bibit[3] L. KOZMA,  On Landsberg  spaces and holonomy of Finsler manifolds,
   \it Contemporary Mathematics \bf 196 \rm(1996), 177-185.

\bibit[4]  G. S.  ASANOV,  Finsleroid--Finsler  space with Berwald and  Landsberg conditions,
  arXiv:math.DG/0603472, 2006.

\bibit[5]  G. S.  ASANOV,   Finsleroid-Finsler  space and spray   coefficients,   arXiv:math.DG/0604526, 2006.

\bibit[6] G. S.  ASANOV, Finsleroid--Finsler spaces of positive--definite and relativistic types,
\it Rep. Math. Phys.  \bf 58 \rm(2006), 275--300.

\bibit[7] G. S. ASANOV,  Finsleroid--Finsler space and geodesic spray    coefficients,
{\it Publ.  Math. Debrecen } {\bf 71/3-4} (2007), 397-412.

\bibit[8] R. S. INGARDEN and L. TAMASSY,
 On parabolic geometry and irreversible macroscopic time,
  \it Rep. Math. Phys. \bf 32 \rm(1993), 11.

\bibit[9] R. S.  INGARDEN,  On physical applications of Finsler geometry,
   \it Contemporary Mathematics \bf 196 \rm(1996), 213--223.

\bibit[10]  G. S.  ASANOV,    Finsleroid  corrects     pressure  and energy   of
  universe. Respective  cosmological equations,
 arXiv:math-ph/0707.3305v1, 2007.

\bigskip

\end{document}